\documentclass[preprint,12pt]{elsarticle}
\usepackage{booktabs}
\usepackage{mathrsfs} 
\usepackage{float}
\usepackage{subcaption}
\usepackage{algorithm}
\usepackage{algorithmic}
\usepackage{pifont}
\usepackage{amsmath}
\usepackage{amsfonts}
\newcommand{\cmark}{\ding{51}}%
\newcommand{\xmark}{\ding{55}}%

\newlength{\maxwidth}
\newcommand{\algalign}[2]
{\makebox[\maxwidth][r]{$#1{}$}${}#2$}

\usepackage{amssymb}
\usepackage{amsmath}
\usepackage{setspace}
\usepackage[symbol]{footmisc}

\journal{Journal}
\begin{document}

\newpage
\begin{frontmatter}
\title{AI-augmented stabilized finite element method}
\author{Sangeeta Yadav\footnote[1]{Corresponding author}}

\author{Sashikumaar Ganesan}
\affiliation[]{{Department of Computational and Data Sciences,},
            addressline={Indian Institute of Science}, 
            city={Bengaluru},
            postcode={560012}, 
            state={Karnataka},
            country={India}}
\begin{abstract}

An artificial intelligence-augmented Streamline Upwind/Petrov-Galerkin finite element scheme (AiStab-FEM) is proposed for solving singularly perturbed partial differential equations. 
In particular, an artificial neural network framework is proposed to predict optimal values for the stabilization parameter. The neural network is trained by minimizing a physics-informed cost function, where the equation's mesh and physical parameters are used as input features. Further, the predicted stabilization parameter is normalized with the gradient of the Galerkin solution to treat the boundary/interior layer region adequately. The proposed approach suppresses the undershoots and overshoots in the stabilized finite element solution and outperforms the existing neural network-based partial differential equation solvers such as Physics-Informed Neural Networks and Variational Neural Networks. 
\end{abstract}

\end{frontmatter}

\section{Introduction}
\label{sec:intro}
 The convection-diffusion partial differential equation (PDE) is used in many applications to model phenomena involving the transport of particles and energy. When the convection dominates the diffusion, the PDE becomes Singularly Perturbed Partial Differential Equation~(SPPDE). The accurate numerical solution of SPPDEs is challenging to find due to the presence of interior and boundary layers. To suppress spurious oscillations in the numerical solution of SPPDEs, standard numerical methods such as finite difference, finite volume, and  finite element method~(FEM) are stabilised with a stabilization method. Several stabilizations and shock-capturing techniques such as slope limiting, flux limiting, artificial dissipation, etc.~\cite{Leveque90, Roos2008} have been proposed in the literature. Adding artificial diffusion is the key idea in almost all finite element stabilization methods, and in particular, the Streamline Upwind Petrov Galerkin (SUPG) stabilization method is most widely used for solving SPPDEs~\cite{Brooks1982, Burman2010,ganesan_operator-splitting_2012, VJohn2011, Behzadi2016, DeFrutos2014, Johnson1984, Giere2015, Li2019}. Although SUPG suppresses spurious oscillations, it requires a well-tuned stabilization parameter to suppress undershoots and overshoots in the solution. Moreover, finding an optimal value of the stabilization parameter is very challenging~\cite{Knobloch2019,knobloch_importance}, and a heuristic approach is often followed.

Several neural network-based solvers, such as Physics-Informed Neural Networks (PINNs)~\cite{raissi2019physics}, have recently been proposed in the literature to solve PDEs without involving a mesh-based numerical method. This technique uses an artificial neural network~(ANN) to estimate the pointwise solution by minimizing the residual of the given differential equation. These neural network methods follow an unsupervised or semi-supervised approach, i.e. labelled data is not used to train the neural network. 
Instead, the residual of the governing PDE is used as the cost function of the neural network to minimize it implicitly. 
Alternatively, the variational-PINNs (VPINNs)~\cite{kharazmi2019variational} relies on  the 
residual of the variational form of the PDEs as the cost function requires relatively lower-order derivatives, resulting in improved accuracy and fewer smoothness conditions on the solution. 
A similar approach known as VarNet~\cite{VarNet} additionally delivers a Reduced-Order-Model (ROM) for efficient computations. Recently, a novel framework called hp-variational PINNs (hp-VPINNs) has been proposed, which allows hp-refinement of the approximate solution~\cite{kharazmi2020hp}. 
Although the neural network-based solvers are more robust than the mesh-based solvers and certainly show promise, it is observed that the  mesh-based methods are efficient and often achieve better accuracy with finer mesh. Further, network-based solvers have limited performance on SPPDEs, and a stabilized form is still needed. This forms the main motivation of the present study, where we augment the SUPG-FEM with deep learning algorithms to obtain an optimal stabilization parameter.
A few ANN-augmented stabilization schemes have also been proposed in the literature. 
Schwander et al.~\cite{Schwander2021} has proposed stabilization of Fourier Spectral Methods (FSMs) with artificial dissipation. The authors use ANN to estimate the regularity of the local solution to identify regions where artificial dissipation can be added. 
Discacciati et al.~\cite{Discacciati2020} used a similar approach for stabilizing Runge-Kutta discontinuous Galerkin (RKDG) methods. 
Veiga and Abgrall~\cite{Veiga2020} have proposed a method for parameter-free stabilization of FEM using ANNs. 
Ray and Hesthaven~\cite{Ray2018, Ray2019} used multilayer perceptrons to identify troubled cells to perform slope limiting for high-order FEM to control spurious oscillations in the solution.
We refer to~\cite{Lye2020} for a detailed review of the ANN-augmented methods in computational fluid dynamics. Moreover, the existing ANN-augmented stabilization methods mostly follow a supervised learning approach, where a reasonable amount of labelled data is needed to train the ANN.

In our earlier study on ANN-augmented stabilization~\cite{yadav2021spde}, we used $L^2$-Error minimization for training an ANN, SPDE-Net. In the $L^2$-Error minimization approach, the network is trained by minimizing the error in the  numerical solution, which requires the analytical solution.  

A similar $L^2$-Error minimization approach is used by Tomasso et al.~\cite{Tassi_2022} for two-dimensional SPPDEs. 
Researchers have developed similar schemes for various numerical methods, such as Fourier Spectral and discontinuous Galerkin.
For example, Lukas et al.~\cite{hesthaven_spectral} proposed a local ANN to estimate the local solution regularity for controlling oscillations and demonstrated the efficiency of nonlinear artificial viscosity methods for Fourier Spectral methods. In another work, Jian et al.~\cite{YU2022105592} showed the ANN's potential for shock capturing in discontinuous Galerkin methods. Their proposed ANN model maps the element-wise solution to a smoothness indicator to determine the artificial viscosity.

In this paper, we propose an ANN-augmented stabilization strategy for SUPG-FEM. In particular, we strive to train the ANN using the residual of the differential equation, which eliminates the need for an analytical solution.  Unlike the existing ANN-augmented stabilized methods, the proposed approach follows unsupervised neural network training. 
The main contributions of this study are summarized below.


 \begin{itemize}
 \item A hybrid approach combining the strengths of the SUPG-FEM and DNN, where a DNN predicts an optimal value of the stabilization parameter and is consequently used in the SUPG-FEM framework to solve SPPDEs.

\item A physics-informed cost function based on aposteriori error estimator~\cite{verfurth} to train the neural network. Unlike the existing approaches in the literature, this cost function does not require the analytical solution of the differential equation.

\item The gradient normalization is applied to the predicted global stabilization parameter ($\hat{\tau}$) to incorporate the dynamics from the interior/boundary layers regions.
\end{itemize}

\section{Preliminaries}\label{sec:preliminaries}
Let $\Omega \subset \mathbb{R}^2$ be a domain bounded with a polygonal Lipschitz-continuous Dirichlet boundary, $\Gamma_\text{D}:=\partial\Omega$. We use the standard notations $\text{L}^\text{p}(\Omega)$ and $\text{W}^{\text{k,p}}(\Omega)$, $1\leq \text{p} <\infty$, $\text{k}\geq 0$, to denote the Lebesgue and Sobolev spaces respectively. Further,
H$^\text{k}(\Omega)$ (i.e. W$^{\text{k,2}}(\Omega)$) denotes the Hilbert space,  $(\cdot,\cdot)$  denotes the inner product in the $\text{L}^2(\Omega)$ space and $|c|$ stands for the Euclidean norm of $c \in \mathbb{R}$.

\subsection{Convection-Diffusion Equation}
 A two-dimensional convection-diffusion equation is given as follows:
\begin{align}\label{eq:convection diffusion}
\begin{split}
-\varepsilon \Delta u + \mathbf{b}\cdot\nabla u &= f \text{ in } \Omega = [0,1]^2,\\ u &= u_b \text{ on } \Gamma_\text{D}, \\
\end{split}
\end{align}
where $\varepsilon > 0 $ is the diffusion coefficient, $\mathbf{b}=(b_1,b_2)^\text{T}$ is the convective velocity, $f \in \text{L}^2(\Omega)$ is an external source term, $u$ is the unknown scalar term. Here, $u_b \in \text{H}^{1/2}(\Gamma_\text{D})$ is a known function.

\subsection{Weak Formulation}
Let $U:= \text{H}^1(\Omega)$ be the solution space.
  We derive weak form of the equation \eqref{eq:convection diffusion} by multiplying it by a test function $v \in V:= \text{H}^1_0(\Omega)$, where ~{$\text{H}^1_0(\Omega) :=\{ v\in U \text{ and } v=0 \text{ on } \Gamma_D\}$},
 followed by integrating on $\Omega$ and subsequently applying integration by parts. The obtained weak form reads as 
 follows:\\
 
  \noindent Find $u\in U$ such that for all $v \in V$
 \begin{align}\label{eq:3}
 a(u,v)  &= (f,v), 
 \end{align}
 where the bilinear form $a(\cdot,\cdot): U \times V \rightarrow \mathbb{R}$ is defined by 
 \begin{align}\label{eq:4}
 a(u,v) &= \int _\Omega \varepsilon \nabla u\cdot \nabla v\,dx +\int_\Omega \textbf{b}\cdot\nabla u v\,dx\\
 (f,v) & = \int _\Omega fv\,dx.
  \end{align}

Let $\Omega_h$ be an admissible decomposition of $\Omega$ and let $K$ be a single
cell in $\Omega_h$. Further, $U_h\subset U$ and  $V_h \subset V$  be finite-dimensional subspaces.
The discrete form of the equation~\eqref{eq:3} reads:\\

\noindent Find $u_h\in U_h$ such that
\begin{equation}\label{eq:9}
\begin{split}
a_h(u_h,v) :=\varepsilon\left(\nabla u_{h}, \nabla v\right)+\left(\mathbf{b} \cdot \nabla u_{h}, v\right)=\left(f, v\right)
\end{split}
\end{equation}
for all $v \in V_h$.
It is well-known that the standard Galerkin form \eqref{eq:9} induces spurious oscillations
in the solution, particularly in the presence of the boundary and internal layers. To avoid these spurious oscillations, the SUPG stabilization term is introduced in the discrete formulation \eqref{eq:9}.

\subsection{SUPG stabilization}
In SUPG, we add a residual term to the variational form in the streamlined direction. 
Let $R(u)$ be the residual of the equation \eqref{eq:convection diffusion} defined as:
\begin{equation}\label{eq:stdtau}
\begin{split}
R(u)  &= -\varepsilon \Delta  u_h + b\cdot \nabla u_h - f.\\
\end{split}
\end{equation}
Now, the SUPG weak form reads as follows. 
\\

\noindent Find $ u_h \in U_h$ such that:
\begin{equation}\label{eq:SUPG_weak_form}
\begin{split}
a^{\text{SUPG}}_h(u_h,v) &= \varepsilon (\nabla u_h, \nabla v ) + (\mathbf{b}\cdot \nabla u_h,v) \\&+ 
 \sum_{K \in \Omega_h}\tau_{K} (-\varepsilon \Delta u_h + \mathbf{b}\cdot \nabla u_h - f_h,\mathbf{b}\cdot \nabla v)_{K} \\
 &= f_h(v)  \quad \forall v \in V_h.\\
 \end{split}
 \end{equation} 
 Here, $\tau_K \in \mathbb{R}$ is a user-chosen non-negative stabilization parameter, and it plays a key role in the accuracy of the approximated solution. A large value results in unexpected smearing, whereas a small value would not suppress the oscillations. Therefore, an optimal stabilization parameter is required to control oscillations and smearing properly. 
\begin{figure*}[h]
\begin{center}
\includegraphics[width=\textwidth]{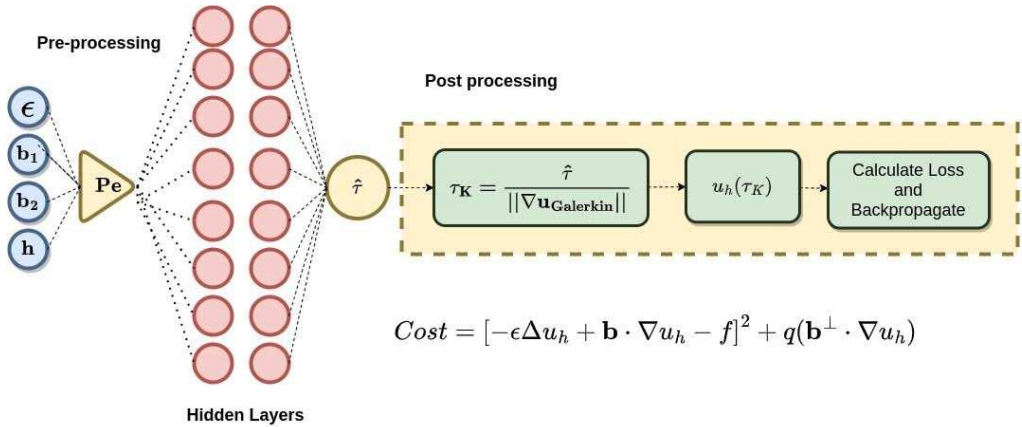}
\caption{AiStab-FEM: Proposed Network Architecture}
\label{fig:net_arch}
\end{center}
\end{figure*}
\section{AiStab-FEM: Optimal stabilization parameter prediction}
The SUPG scheme suffers from a major limitation, viz., the optimal global stabilization parameter $\tau$
is not easy to find for a general case. 
While a standard expression is available for the scalar convection-diffusion equation, a general expression is still unknown. This necessitates devising a technique that adaptively finds an optimal value of the stabilization parameter. 
Research in machine learning, mainly ANNs, has seen massive growth in the last few years. 
ANNs have been used in virtually every scientific field due to their `universal approximator' nature~\cite{Hornik1989}.
In this article, we exploit this property of ANNs to approximate the SUPG stabilization parameter. 
We propose a neural network-based novel framework, AiStab-FEM, for finding a near-optimal stabilization parameter. 
\subsection{AiStab-FEM: $\tau$ approximation}
 The problem of approximating the stabilization parameter~($\tau$) is modelled as a regression problem and solved using a neural network called AiStab-FEM.  We explain AiStab-FEM coherently in algorithm~\ref{alg:AiStab-FEM}.

\begin{algorithm}[H]
\setstretch{2.07}
\caption{AiStab-FEM}
\begin{algorithmic}[1]
\STATE{Given: \{I: input, $\eta$: learning rate, $n_{epochs}$: number of training epochs}\}
\STATE{Initialize the neural network weights with random parameters $\theta_0$}
\STATE{Initialize the optimizer (Adam in this case) and stepLR scheduler}
\STATE{I = $\{\varepsilon, b_1,b_2,h\}$}
\FOR {$t = 1$ to $n_{epochs}$} 
\settowidth{\maxwidth}{$m$}
\STATE {$\displaystyle Pe(I)  =  \frac{|\mathbf{b}|h}{2\varepsilon}$}
\STATE{$\hat{\tau}(\theta_t) = \text{AiStab-FEM}_{\theta_t}(Pe(I))$}
 \STATE{$\displaystyle \tau_K = \frac{\hat{\tau}(\theta_t)}{||\nabla u_{Galerkin}(I)||}$}
 \STATE{$q(s) = 
\begin{cases}
\sqrt{s} & s>1  \\
 2.5 s^2 -1.5s^3 &\text{otherwise}
\end{cases}$}
\STATE{ $\mathbf{b}^\perp = 
\begin{cases}
\displaystyle \frac{(b_2,-b_1)}{|\mathbf{b}|} & \text{when } |\mathbf{b}|\neq 0 \\
0 & \text{when } |\mathbf{b}| = 0
\end{cases}$}
\STATE{Solve equation \eqref{eq:SUPG_weak_form} with $\tau_K$ and get $u_h$}
\STATE{$\displaystyle cost({\theta_t}) = \sum_{K \in \Omega} ||-\varepsilon \Delta u_h + \mathbf{b}\cdot \nabla u_h -  f||_{K}^2
+~q(\mathbf{b}^\perp \cdot \nabla u_h)$}  
\STATE{$\displaystyle \theta_{t+1} = \theta_t - \eta_t  \nabla cost(\theta_t)$}
\ENDFOR
\end{algorithmic}
\label{alg:AiStab-FEM}
\end{algorithm}
 The input vector I, to the AiStab-FEM scheme, comprises of diffusion coefficient ($\varepsilon$), convection velocity ($\mathbf{b}$) 
 and mesh size ($h$). In algorithm~\ref{alg:AiStab-FEM}, we calculate P\'eclet (Pe) number and the Galerkin solution for the input sample I. 
A global $\hat{\tau}$ for all the cells is inadequate since layered regions should get a different stabilization value than the other cells. To address this issue, we normalize the predicted $\hat{\tau}$ by dividing it by the cell-wise Euclidean norm of the gradient of the solution calculated with the Galerkin technique. It is to be noted that this Galerkin solution tends to be sub-optimal for a general case. In algorithm \eqref{alg:AiStab-FEM}, $u_h$ is the SUPG solution of equation \eqref{eq:SUPG_weak_form}, which is used to calculate the cost function. 

\subsection{Cost formulation}
The cost function is a very important component of neural network training. It should be an appropriate metric for the error in the prediction. We consider the residual-based error indicator proposed by~\cite{JOHN20112916},~\cite{Knobloch2019} as the cost function. It consists of a residual-based term and cross-wind derivative, as shown below: 
\begin{equation}\label{eq:loss}
\setlength{\jot}{8pt}
\begin{split}
cost(u_h(\theta)) &= \left[-\varepsilon \Delta u_h+\mathbf{b}\cdot \nabla u_h -  f\right]^2 + q(\mathbf{b}^\perp \cdot \nabla u_h)\\
\theta^* &= \arg \min_{\theta}\left( cost(u_h(\theta))\right) \\
\end{split}
\end{equation}
$q(s)$ and $\mathbf{b}^\perp$ are defined in Algorithm \eqref{alg:AiStab-FEM}, AiStab-FEM aims to find $\theta^*$ so that the numerical error calculated with predicted $\hat{\tau}$ is minimal.
We find optimal $\theta^*$ by back-propagating this cost (equation \eqref{eq:loss}) in an iterative manner.

\subsection{Network Architecture}
We show the architecture of the proposed network AiStab-FEM diagrammatically in figure \ref{fig:net_arch}. It is a fully connected neural network with two hidden layers with an input layer of size four. We use sigmoid non-linearity for the output layer and $\tanh$ after each hidden layer. We use the Adam optimizer to accelerate the gradient descent algorithm and the StepLR scheduler to control the learning rate adaptively. It decays the learning rate to gamma times the original learning rate at every step. Other hyperparameters for training are given in~Table \ref{tab:hyperparameters}.  
\vskip 2pt
\begin{table}[!h]\centering
\vskip -10pt
\caption{Network hyper-parameters}\label{tab:hyperparameters}
\begin{tabular}{lrr}\toprule
Hyper-parameter&Value\\
\midrule
Non-linearity &tanh \\
Optimizer &Adam \\
Initial LR &0.0001 \\
Scheduler &LR \\
Gamma &0.1 \\
Size of the hidden layers &[16,16] \\
\bottomrule
\end{tabular}
\end{table}

\subsection{Baselines}\label{sec:baselines}
To benchmark the proposed technique, we compare it with a few baseline techniques, as explained in this section.  We consider two purely neural network-based PDE solvers and two techniques for approximating $\tau$, which can be used in a mesh-based solver. This way, we ensure a fair comparison of the proposed technique with the pure neural network-based solvers and the mesh-based stabilized solvers. The baseline techniques are as follows:


\begin{enumerate}
\item \textbf{PINNs:} This is a neural network trained for solving PDEs by minimizing the residual of the given equation. For further details, refer~\cite{raissi2017physicsI}.
\item \textbf{Variational Neural Networks for the Solution of Partial Differential Equations~(VarNet):}
In this technique, the numerical solution of the equation is approximated by minimizing the residual of the weak form of the equation. For further details, refer to~\cite{VarNet}.
\item \textbf{Standard $\tau_{std}$:} Given $Pe$ is the local P\'eclet number as given in Algorithm \ref{alg:AiStab-FEM}, we use the following expression to compute $\tau_{std}$:
\begin{equation}\label{eq:std_tau}
\tau_{std} =\frac{h}{2|\mathbf{b}|}\left(\coth{(Pe)} - \frac{1}{Pe}\right)\\
\end{equation}
where $\tau_{std}$ is used to solve SUPG weak form \eqref{eq:SUPG_weak_form} of the given equation. 
Table \ref{tab:std_tau_values} shows the value of $\tau_{std}$ from the standard formula for different testing examples. These values will be used for the baseline comparison. 
\item \textbf{Standard $\tau_{std}$ normalized with $||\nabla u_{Galerkin}||$:}
In this technique, we normalize the $\tau_{std}$~(equation \eqref{eq:std_tau}) with the cell-wise Euclidean norm of the gradient of the numerical solution.
\begin{equation}
\begin{split}
\tau_{std} &=\displaystyle \frac{h}{2|\mathbf{b}|} \left(\coth{(Pe)} - \frac{1}{Pe}\right)\\
\text{where }\tau_K &= \displaystyle \frac{\tau_{std}}{||\nabla u_{Galerkin}||}
\end{split}
\end{equation}
As is shown in algorithm \ref{alg:AiStab-FEM}, the proposed technique involves the normalization of $\hat{\tau}$ with $||\nabla u_{Galerkin}||$, so we normalize the $\tau_{std}$ as well to make a fair comparison against the standard technique.
\end{enumerate}

\subsection{Error Metrics}\label{sec:error_metrics}
We use the following metrics to calculate numerical errors in the solution obtained with the predicted $\hat{\tau}$. We will use them for comparison against the baselines as explained in section~[\ref{sec:baselines}]

\begin{equation}
\begin{aligned}
L^2\text{-error: } & ||e_h||_{L^2(\Omega)}  = ||u_h(\hat{\tau})-u||_{L^2(\Omega)} =  \left(\int_{\Omega}(u_h(\hat{\tau})-u)^2dx\right)^{\frac{1}{2}}\\
l^2-\text{error(relative): }\displaystyle & \sum_{i=1}^{N}\frac{||u_h(\hat{\tau})(x_i)-u(x_i)||_0}{||u||_0}  :x_i\in \Omega\\
H^1\text{-error: } & ||e||_{H^1(\Omega)} = \sum_{|\alpha|<=1}(D^\alpha (u_h(\hat{\tau})-u),D^\alpha (u_h(\hat{\tau})-u)) \\
&= \sum_{\alpha<=1}\int_{\Omega} D^\alpha (u_h(\hat{\tau})-u)D^\alpha (u_h(\hat{\tau})-u)dx\\
L^\infty\text{-error: } & ||e||_{L^\infty(\Omega)} = ess\quad sup\{|u_h(\hat{\tau})-u|:x\in \Omega\}.\\
\end{aligned}
\end{equation}
Here, $\hat{\tau}$ is the stabilization parameter predicted by AiStab-FEM, $u$ is the known analytical solution, $u_h(\hat{\tau})$ is the SUPG solution calculated with $\hat{\tau}$, $D^\alpha$ is the weak derivative (refer definition 1.14 from \cite{Ganesanbook}).

\section{Experiments and Results}
This section discusses the experiments conducted to evaluate AiStab-FEM's performance. We consider four singularly perturbed benchmark examples with varying levels of complexity, each one of which is a highly convective case. We consider a unit square domain $\Omega$ uniformly divided into $2*N_{cell}*N_{cell}$ triangular cells, where $N_{cell}$ is the number of cells in the horizontal direction. These cells are non-overlapping and regular in shape. We consider $P_2(\Omega)$ space for approximating $u_h$ and $DG_0(\Omega)$ for the stabilization parameter($\tau_K$)~(for definition of $P_2(\Omega)$ and $DG_0(\Omega)$ finite elements, refer section 3.3 of \cite{Ganesanbook}).

\subsection{Example 1:} 
\begin{figure}[h!]
    \begin{center}
	\begin{subfigure}[t]{0.5 \textwidth}
        \begin{center}{\includegraphics[width = \textwidth]{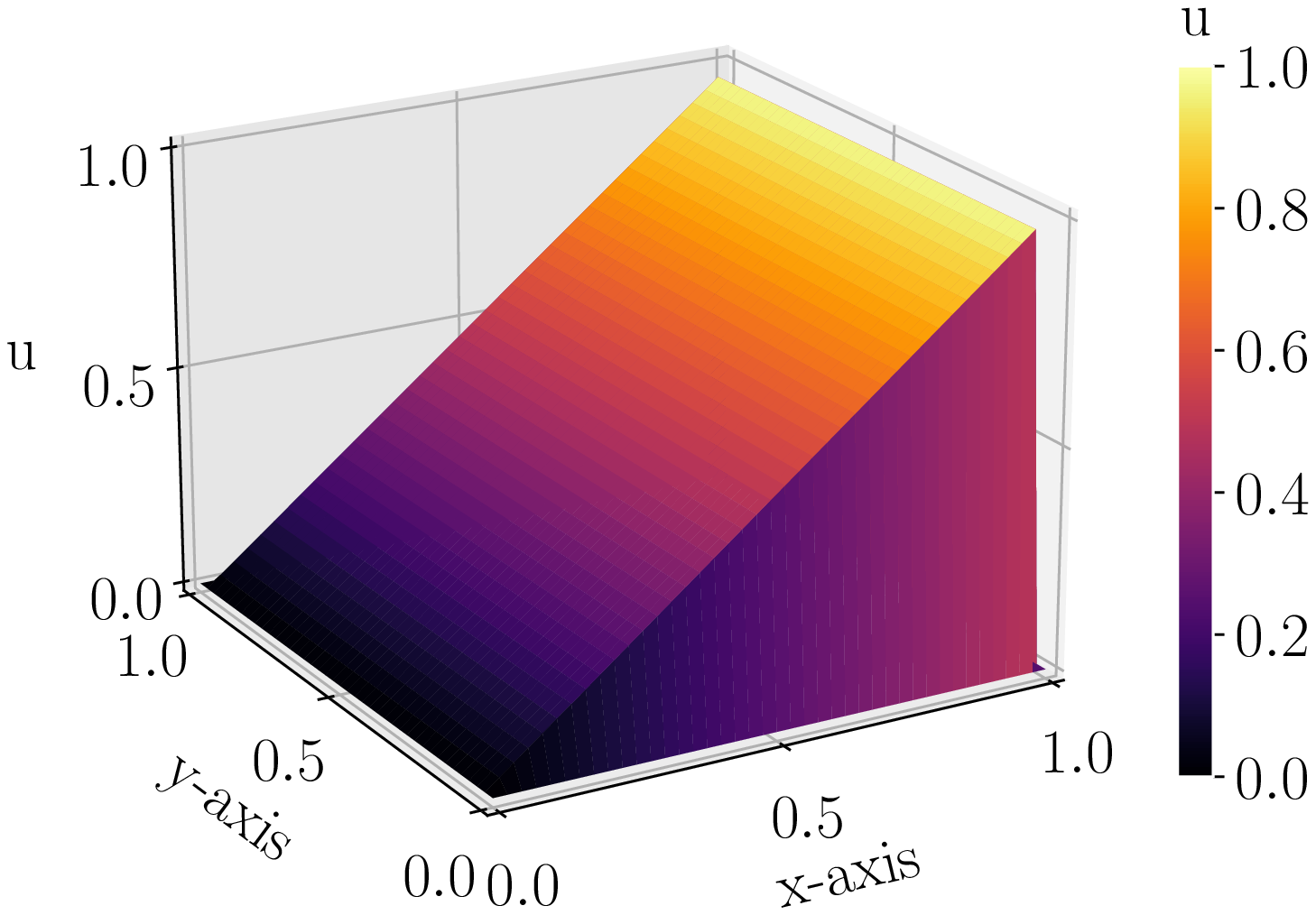}}
	\end{center}
    \caption{\label{Fig:exact1} Exact Solution~($u$)}
\end{subfigure}%
	\begin{subfigure}[t]{0.5 \textwidth}
		\begin{center}{\includegraphics[width=\textwidth]{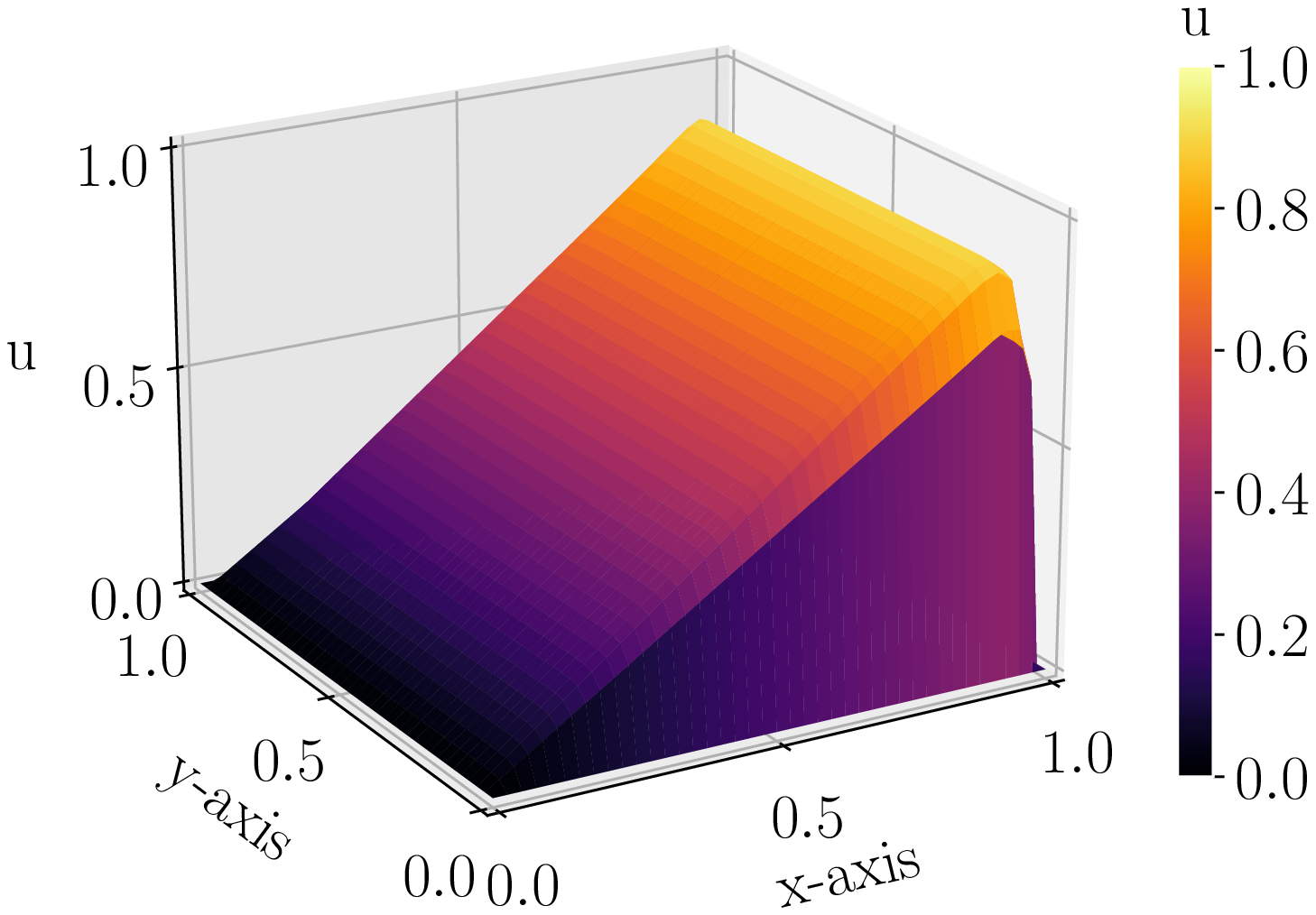}}
		\end{center}
    \caption{\label{Fig:std1} Solution with Standard $\tau_{std}$}
	\end{subfigure}%
	\newline
	\begin{subfigure}[t]{0.5 \textwidth}
		\begin{center}{\includegraphics[width=\textwidth]{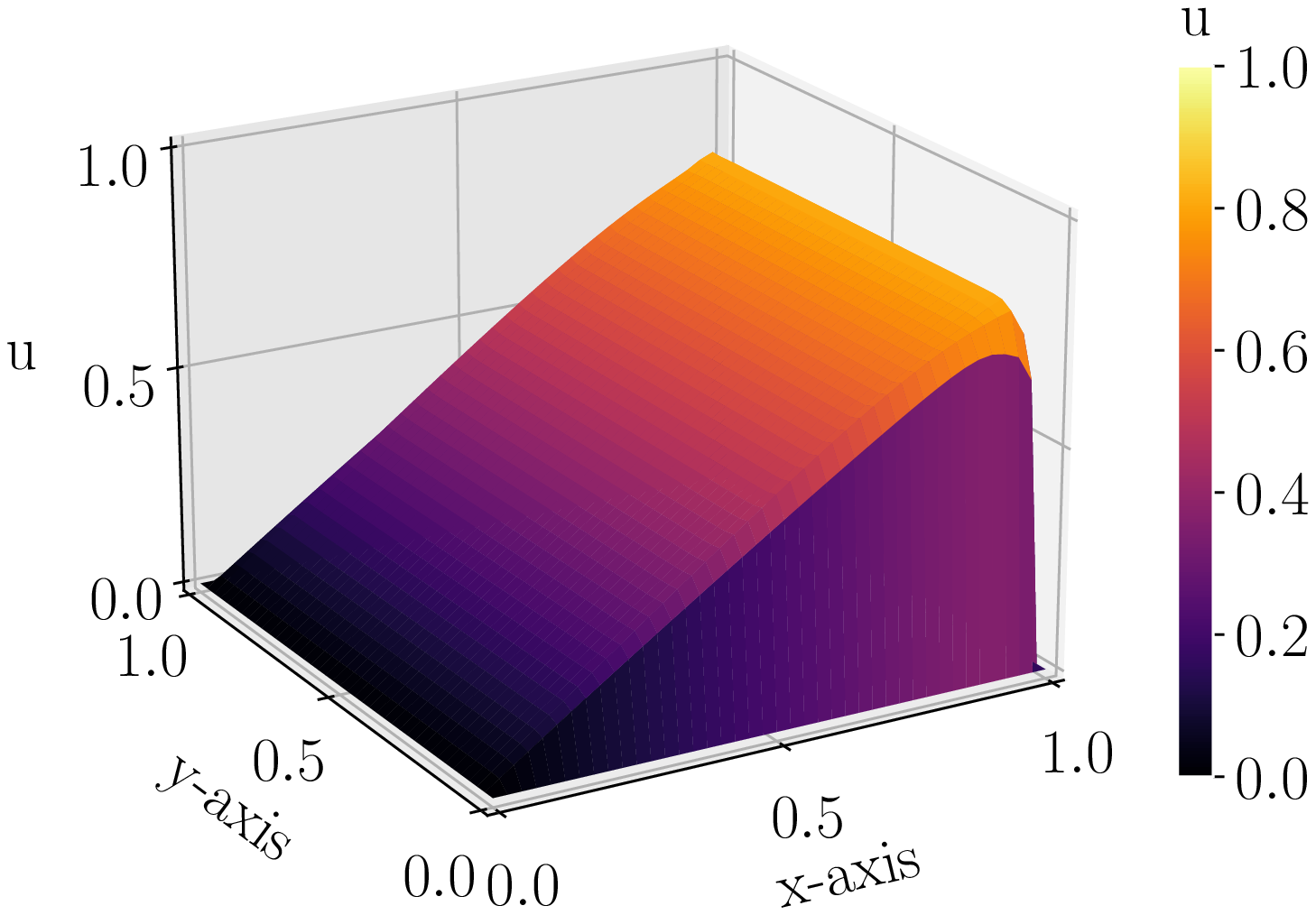}}
		\end{center}
    \caption{\label{Fig:sol1} Solution from AiStab-FEM~($u_h$)}
	\end{subfigure}%
	\begin{subfigure}[t]{0.5 \textwidth}
		\begin{center}{\includegraphics[width=\textwidth]{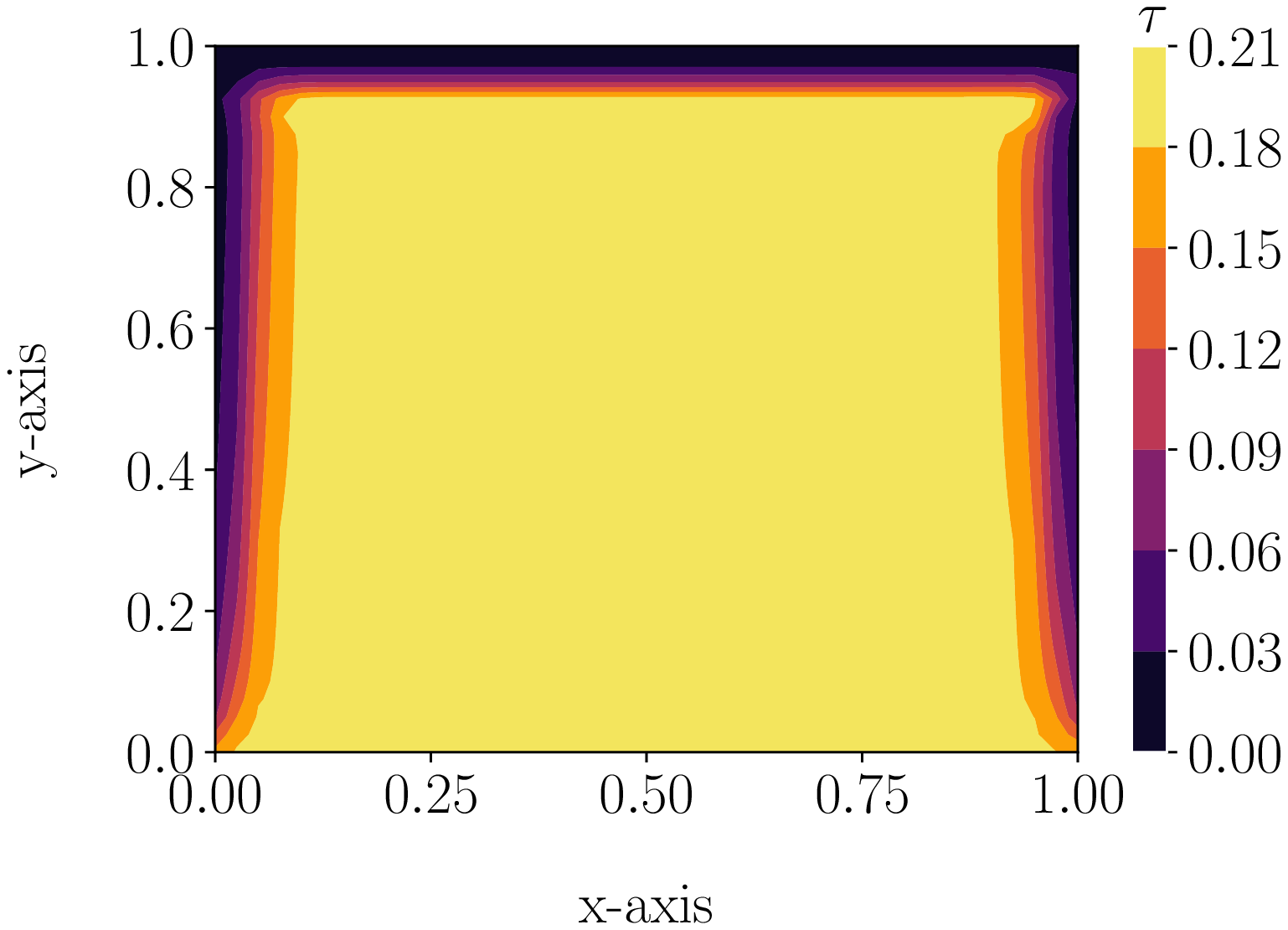}}
		\end{center}
    \caption{\label{Fig:tau1} $\tau_K$ predicted from AiStab-FEM}
	\end{subfigure}%
\end{center}
	\caption{\label{fig:Example1} Example 1 }
\end{figure}

\noindent We consider the equation \eqref{eq:convection diffusion} with following data:
\begin{equation}\label{eq: example1}
\varepsilon = 10^{-8},\quad \mathbf{b}  = (1, 0),\quad
f = 1,\quad \Omega = [0, 1]^2,~u_b = 0.
\end{equation}
This example is taken from~\cite{article2}. 
The exact solution of equation \eqref{eq: example1} is shown graphically in figure \ref{fig:Example1} and is given as follows:
\begin{equation}\label{eq: example1_exact}
\begin{split}
& u(x) = x \quad \text{on} \quad \Omega/\partial \Omega \\
&u(x) = 0\quad \text{on} \quad \Gamma_D\\
\end{split}    
\end{equation}
 $u$ has an exponential layer at $x=1$ and two parabolic layers at $y = 0$ and $y =1$ respectively.
We feed $I = [\varepsilon, b_1,b_2,h]$  to the AiStab-FEM and the predicted $\hat{\tau}$ is normalized with $||\nabla u_{Galerkin}||$. To solve the discretized form of the equation, we consider $40$ cells in the horizontal direction and hence $h_K = \sqrt{2}/40$. We obtain the numerical solution by solving stabilized SUPG weak form of the equation \eqref{eq:SUPG_weak_form} with the predicted $\hat{\tau}$. We use $u$ to compute  error metrics given in section \ref{sec:error_metrics}. 
 The results of AiStab-FEM, for example 1, are shown in figure~\ref{fig:Example1}. In figure \ref{fig:Example1}(d), the value of $\tau_K$, for example 1, is shown on the domain $\Omega$. In the non-boundary region, a constant value is predicted, whereas a gradually decreasing value towards the outer boundary is seen in the boundary region. The highest value of $\tau_K$ is observed at the intersection of two boundary layers. Figure \ref{fig:Example1}(b,~c)  also shows the solution with $\tau_{std}$ normalized by $||\nabla u_{Galerkin}||$ and the solution with AiStab-FEM. Overall, AiStab-FEM outperforms all the baseline techniques (sec. \ref{sec:baselines}) in terms of all the error metrics (sec. \ref{sec:error_metrics}) as shown in  Table \ref{tab:l2error}, \ref{tab:rel2error},
\ref{tab:h1error}, 
\ref{tab:linfinityerror}.

\begin{table}[!htp]\centering
\caption{Value of $\tau$ by standard technique for different examples }\label{tab:std_tau_values}
\begin{tabular}{ccccc}\toprule
Example	& Internal Layer & Boundary Layer & $N_{cells}$	&$\tau_{std}$	\\\midrule
1 &\xmark & \cmark	&40	&1.77e-2	\\
2 &\xmark & \cmark	&40	&1.77e-2	\\
3 & \cmark & \cmark	&40	&1.77e-2	\\
4 & \cmark &\xmark	&40	&4.90e-3	\\
\bottomrule
\end{tabular}
\end{table}

\begin{figure}[h!]
	\begin{center}
	\begin{subfigure}[t]{0.5 \textwidth}
		\begin{center}{\includegraphics[width=\textwidth]{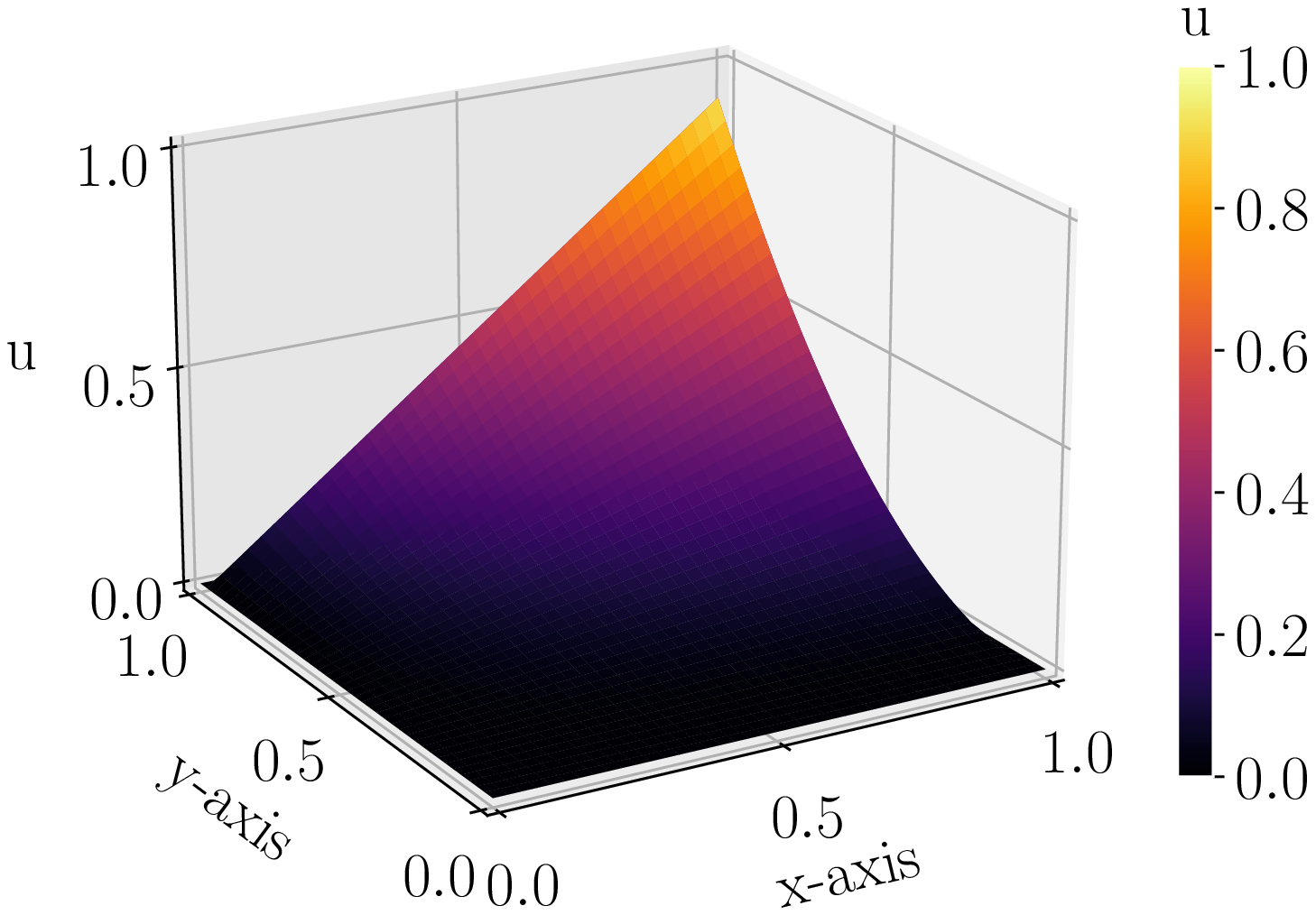}}
		\end{center}
    \caption{\label{Fig:exact2} Exact Solution~($u$)}
	\end{subfigure}%
	\begin{subfigure}[t]{0.5 \textwidth}
		\begin{center}{\includegraphics[width=\textwidth]{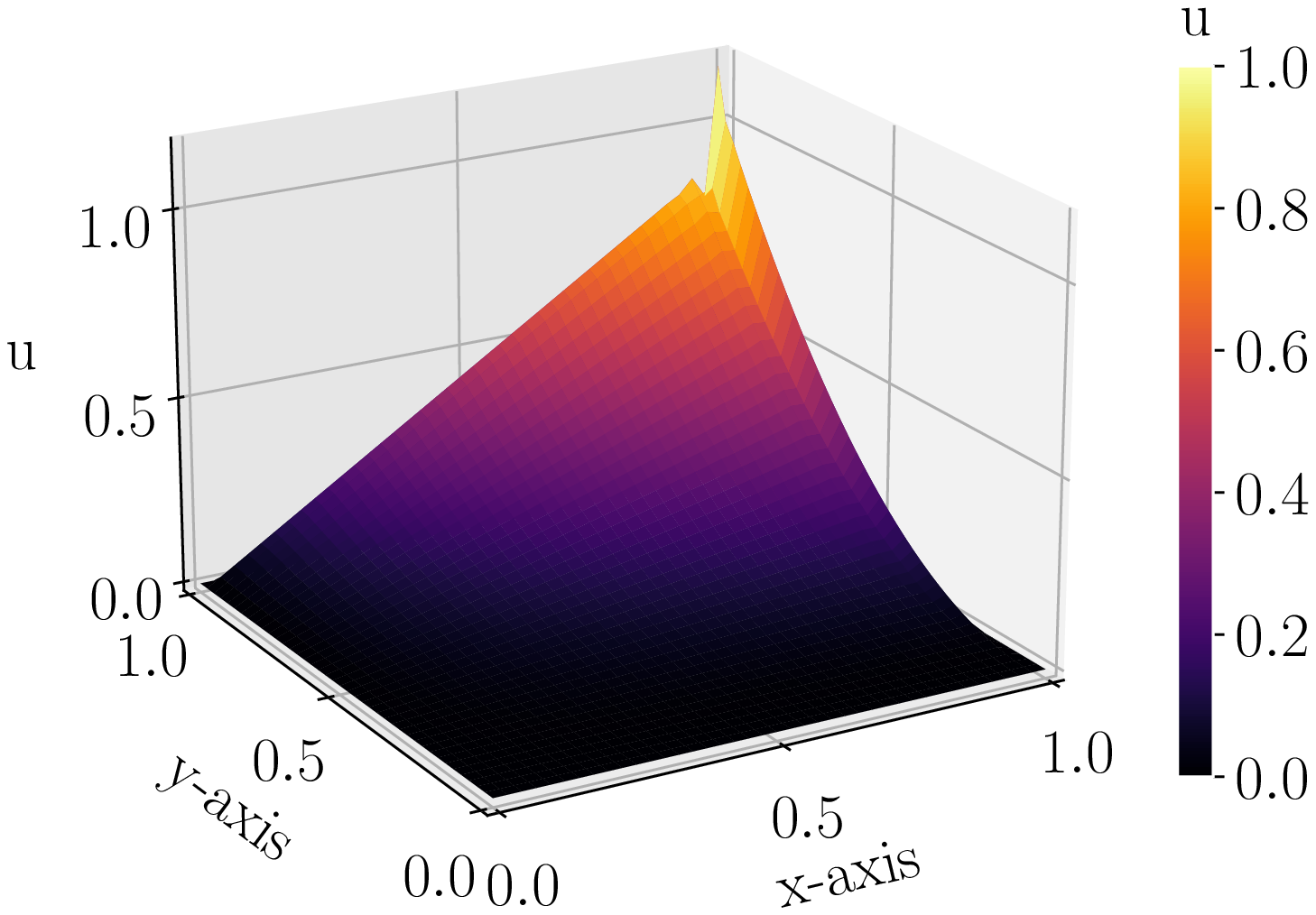}}
		\end{center}
    \caption{\label{Fig:std2} Solution with Standard $\tau_{std}$}
	\end{subfigure}%
	\newline
	\begin{subfigure}[t]{0.5 \textwidth}
\begin{center}{\includegraphics[width=\textwidth]{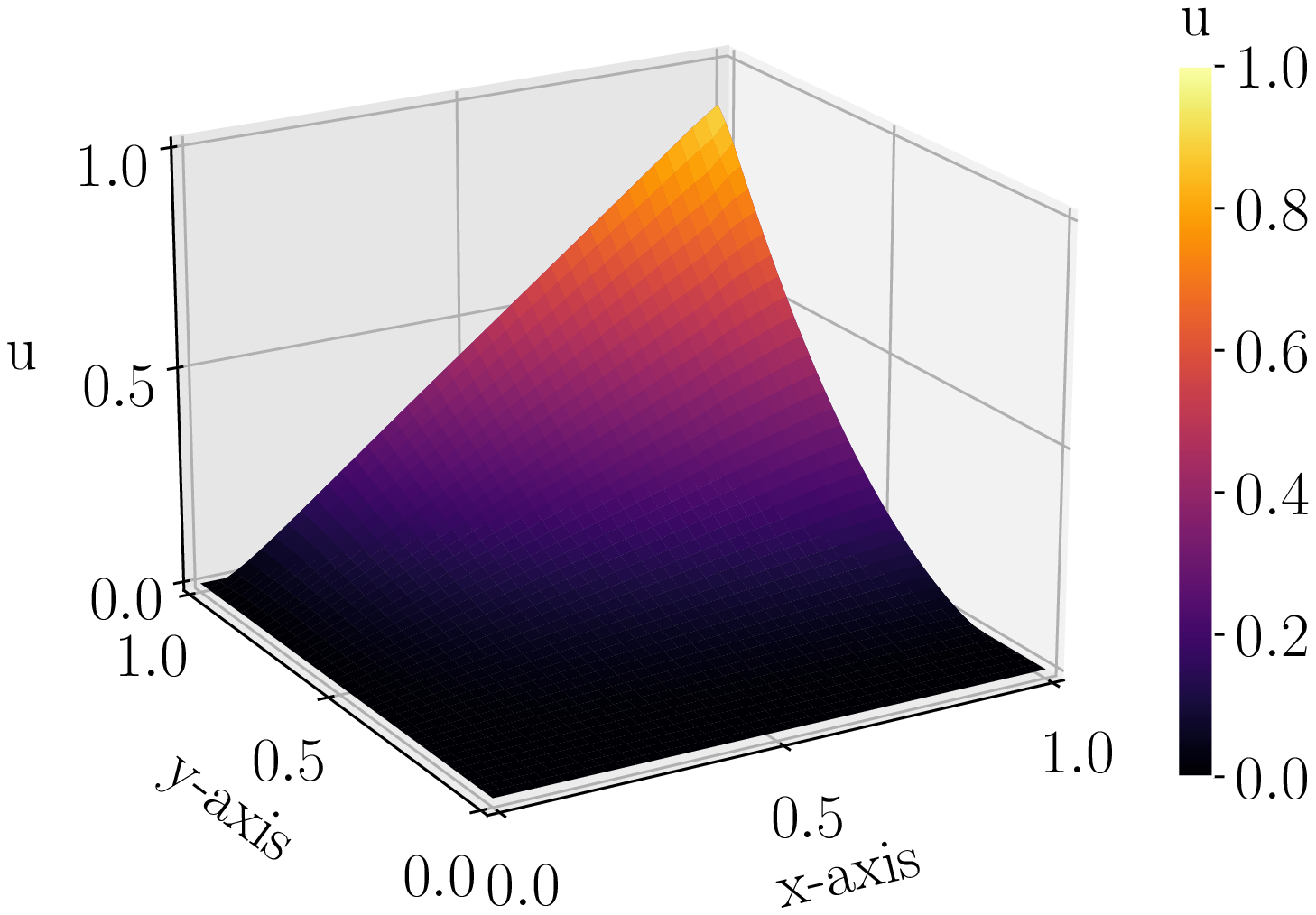}}
\end{center}
    \caption{\label{Fig:sol2} Solution from AiStab-FEM~($u_h$)}
	\end{subfigure}%
	\begin{subfigure}[t]{0.5 \textwidth}
		\begin{center}{\includegraphics[width=\textwidth]{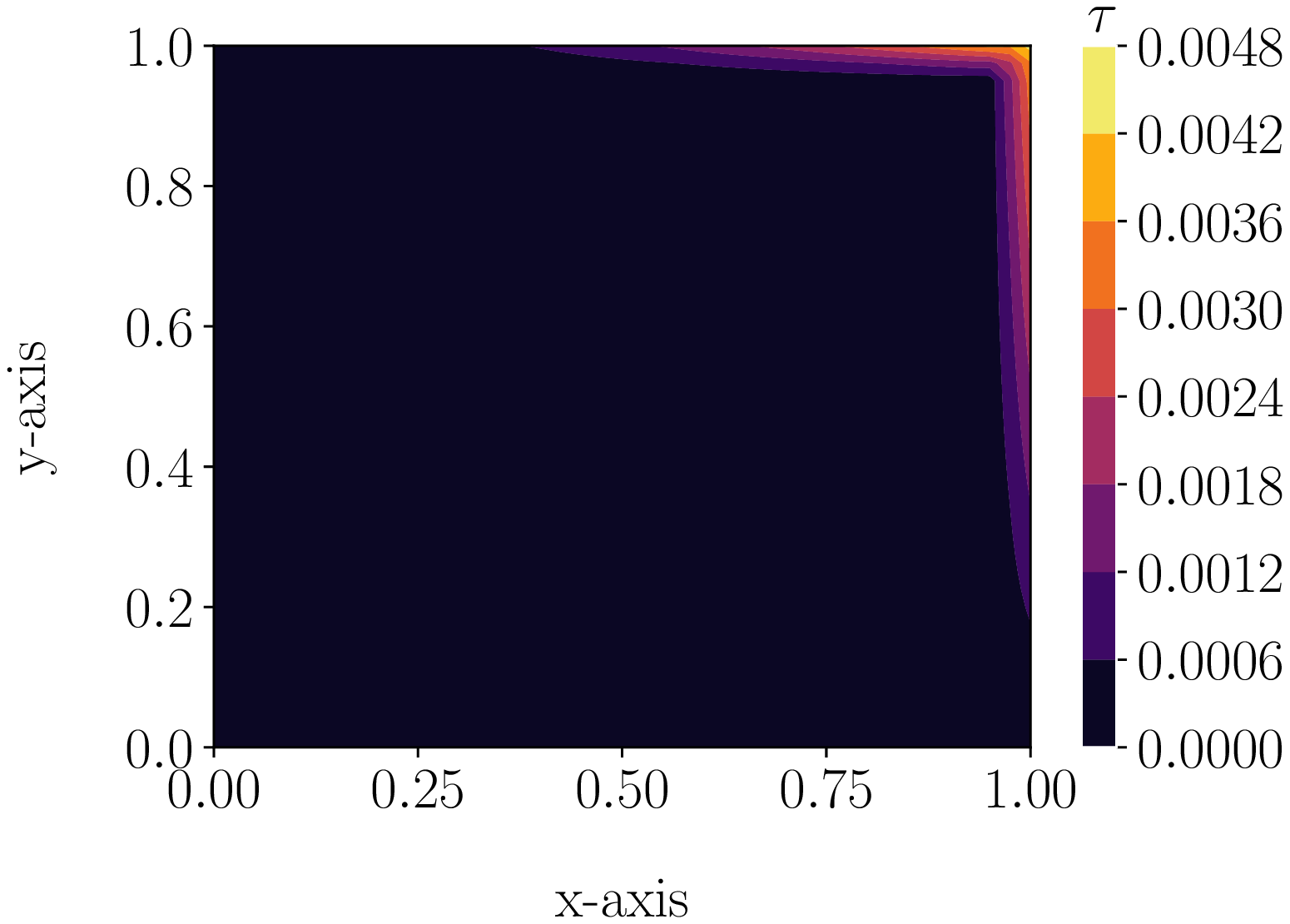}}
		\end{center}
    \caption{\label{Fig:tau2} $\tau_K$ predicted from AiStab-FEM}
	\end{subfigure}%
\end{center}
	\caption{\label{fig:Example2} Example 2 }
\end{figure}

\subsection{Example 2:}
\noindent We consider the convection diffusion equation \eqref{eq:convection diffusion} with following equation coefficients and boundary conditions as given in~\cite{knobloch2009}:
\begin{equation}
\label{eq: example2}
\begin{split}
&\varepsilon = 10^{-8}, \mathbf{b} = (2, 3), u_b = 0\\
\end{split}
\end{equation}
Source term $f$ is calculated by substituting the following analytical solution($u$) in equation \eqref{eq: example2}.
\begin{equation}
\label{eq: example2exact}
\begin{split}
&u(x,y) = x y^{2} - x\exp\left(\frac{3(y - 1)}{\varepsilon} \right) - y^{2}
\exp\left(\frac{2(x - 1)}{\varepsilon} \right) \\
&+ \exp\left(\frac{2(x-1) + 3(y - 1)}{\varepsilon} \right)  
\end{split}    
\end{equation}
It contains two outflow boundary layers near $x = 1.0$ and $y=1.0$ as shown in figure \ref{fig:Example2}(a) and hence makes a suitable test case for checking the performance of the AiStab-FEM as mentioned in~\cite{knobloch2009}. 
The predicted $\tau_K$ is plotted in figure \ref{fig:Example2}(d) and the value of $\tau_{std}$ is $0.0177$ as shown in Table \ref{tab:std_tau_values}. We observe that all the predicted values are smaller than the $\tau_{std}$. We observe a clear pattern of the boundary layers in the heat map of the $\tau_K$. The value of $\tau_K$ is highest around the boundary layer. AiStab-FEM outperforms all the other techniques in terms of all the error metrics for this example, as it does for example 1.

\subsection{Example 3:}
\begin{figure}[h!]
	\begin{center}
	\begin{subfigure}[t]{0.5 \textwidth}
		\begin{center}{\includegraphics[width=\textwidth]{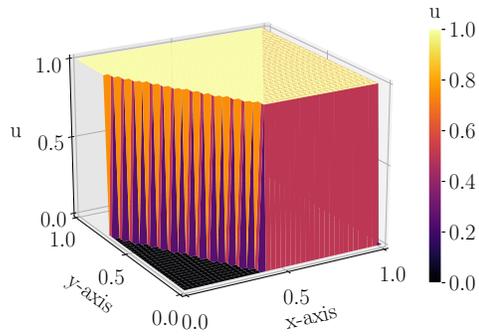}}
		\end{center}
    \caption{\label{Fig:exact3} Exact Solution~($u$)}
	\end{subfigure}%
	\begin{subfigure}[t]{0.5 \textwidth}
		\begin{center}{\includegraphics[width=\textwidth]{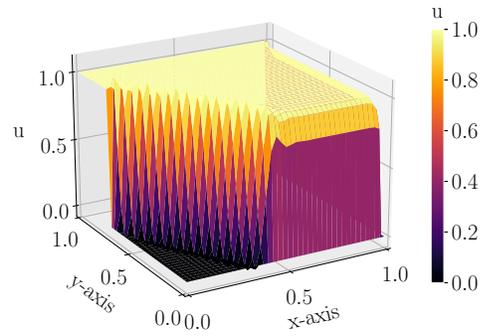}}
		\end{center}
    \caption{\label{Fig:std3} Solution with Standard $\tau_{std}$}
	\end{subfigure}%
	\newline
	\begin{subfigure}[t]{0.5 \textwidth}
		\begin{center}{\includegraphics[width=\textwidth]{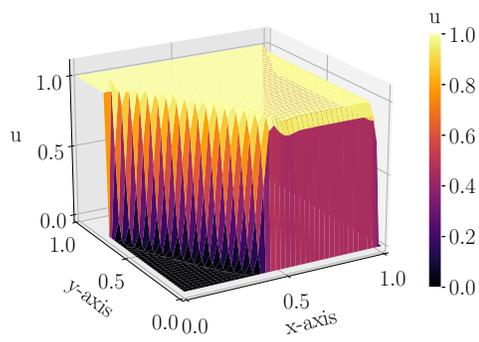}}
		\end{center}
    \caption{\label{Fig:sol3} Solution from AiStab-FEM~($u_h$)}
	\end{subfigure}%
	\begin{subfigure}[t]{0.5 \textwidth}
		\begin{center}{\includegraphics[width=\textwidth]{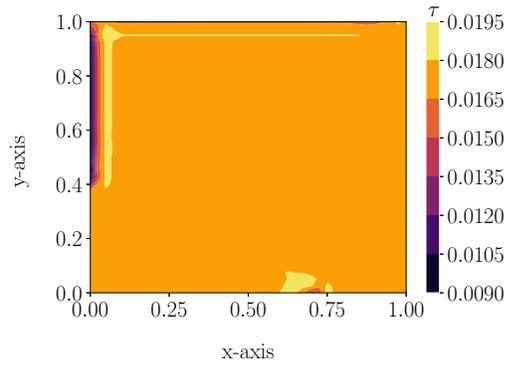}}
		\end{center}
    \caption{\label{Fig:tau3} $\tau_K$ predicted from AiStab-FEM}
	\end{subfigure}%
\end{center}
	\caption{\label{fig:Example3} Example 3 }
\end{figure}

\noindent Next, we consider the convection-diffusion equation \eqref{eq:convection diffusion} with following equation coefficients and boundary conditions:
\begin{equation}
\label{eq: example3}
\begin{split}
&\varepsilon = 10^{-8}, ~\theta = -\pi/3,
~\mathbf{b} = (\cos(\theta),\sin(\theta)),
~f = 0.0,\\
&u_b = \begin{cases}
0,& \text{for } x = 1 \text{ or } y\leq0.7\\
1,              & \text{otherwise}
\end{cases}\\
\end{split}
\end{equation}
This example contained both exponential and boundary layers and was used in~\cite{JOHN20112916}. The value of $\tau_{std}$ for this example is $0.0049$, whereas the value of $\tau_K$ is less than $0.002$ everywhere, as shown in figure \ref{fig:Example3}. We can see that the numerical solution with $\tau_K$ has fewer oscillations than the standard technique. Both the interior and exponential layers are well captured, as can be seen in the heatmap of $\tau_K$ (figure \ref{fig:Example3}(d)). The solutions with the $\tau_{std}$ technique and $\tau_K$ are shown for comparison, and we observe that both result in visible oscillations near the interior layer. We compare  figure \ref{fig:Example3}(b) and \ref{fig:Example3}(c) and observe that a sharp interior layer is captured by AiStab-FEM relative to $\tau_{std}$. AiStab-FEM gives the least errors in terms of all the error metrics~(section \ref{sec:error_metrics}) among all the baseline techniques~(section \ref{sec:baselines}) used for this example.

\newpage
\subsection{Example 4:}
\noindent We consider the convection-diffusion equation \eqref{eq:convection diffusion} with following equation coefficients and boundary conditions:
\begin{equation}
\label{eq: example4}
\begin{split}
&\varepsilon = 10^{-8},~\mathbf{b} = (1,0),~u_b = 0\\
&f = \begin{cases}
0,& \text{ if } |x-0.5|\geq 0.25 \cup |y-0.5|\geq 0.25\\
-32(x-0.5), & \text{otherwise}
\end{cases}\\
&u = \begin{cases}
0,& \text{ if } |x-0.5|\geq 0.25 \cup |y-0.5|\geq 0.25\\
-16(x-0.25)(y-0.75), & \text{otherwise}
\end{cases}
\end{split}
\end{equation}

\clearpage
\begin{figure}[h!]
	\begin{center}
	\begin{subfigure}[t]{0.5 \textwidth}
		\begin{center}{\includegraphics[width=\textwidth]{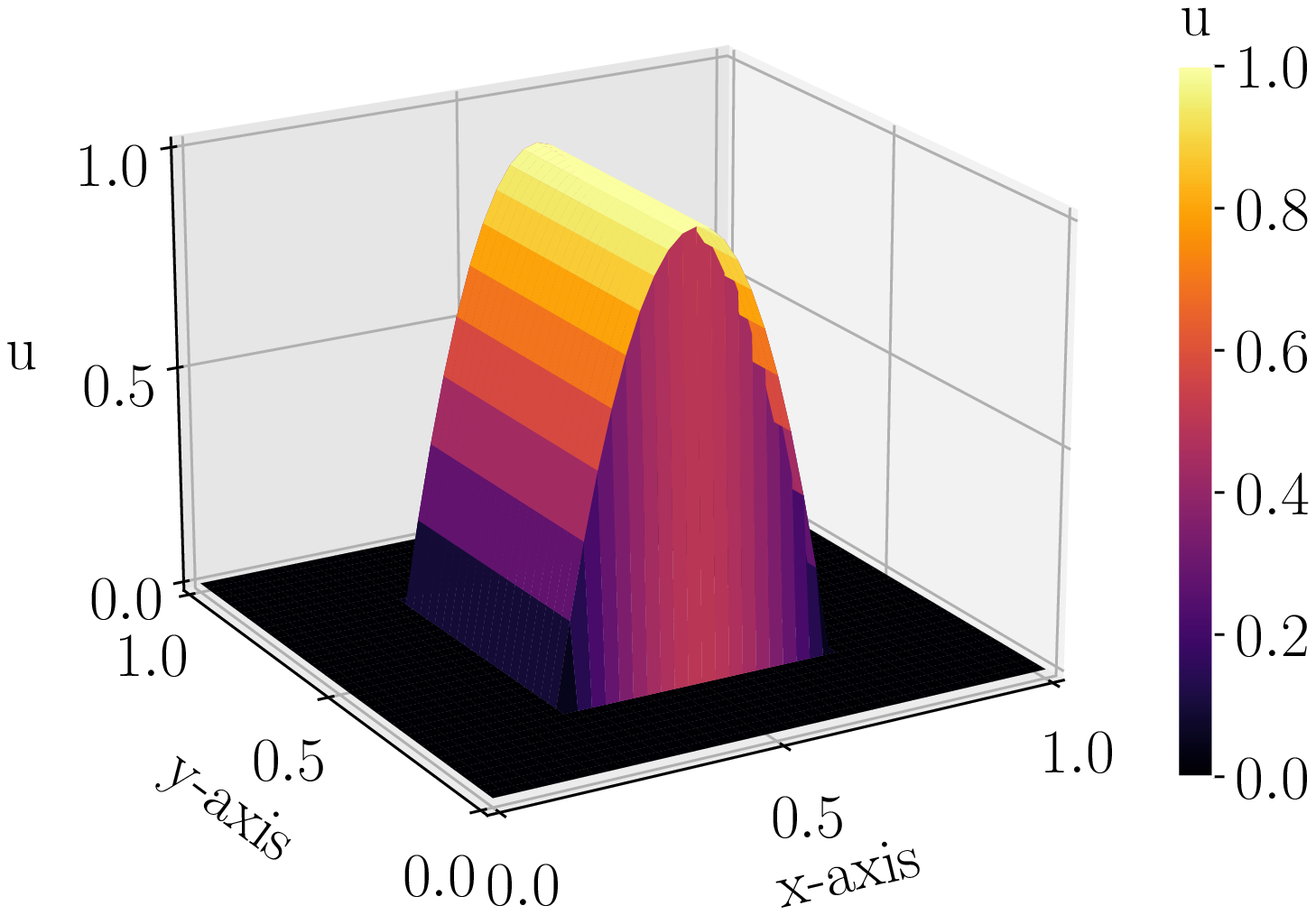}}
		\end{center}
    \caption{\label{Fig:exact4} Exact Solution~(u)}
	\end{subfigure}%
	\begin{subfigure}[t]{0.5 \textwidth}
		\begin{center}{\includegraphics[width=\textwidth]{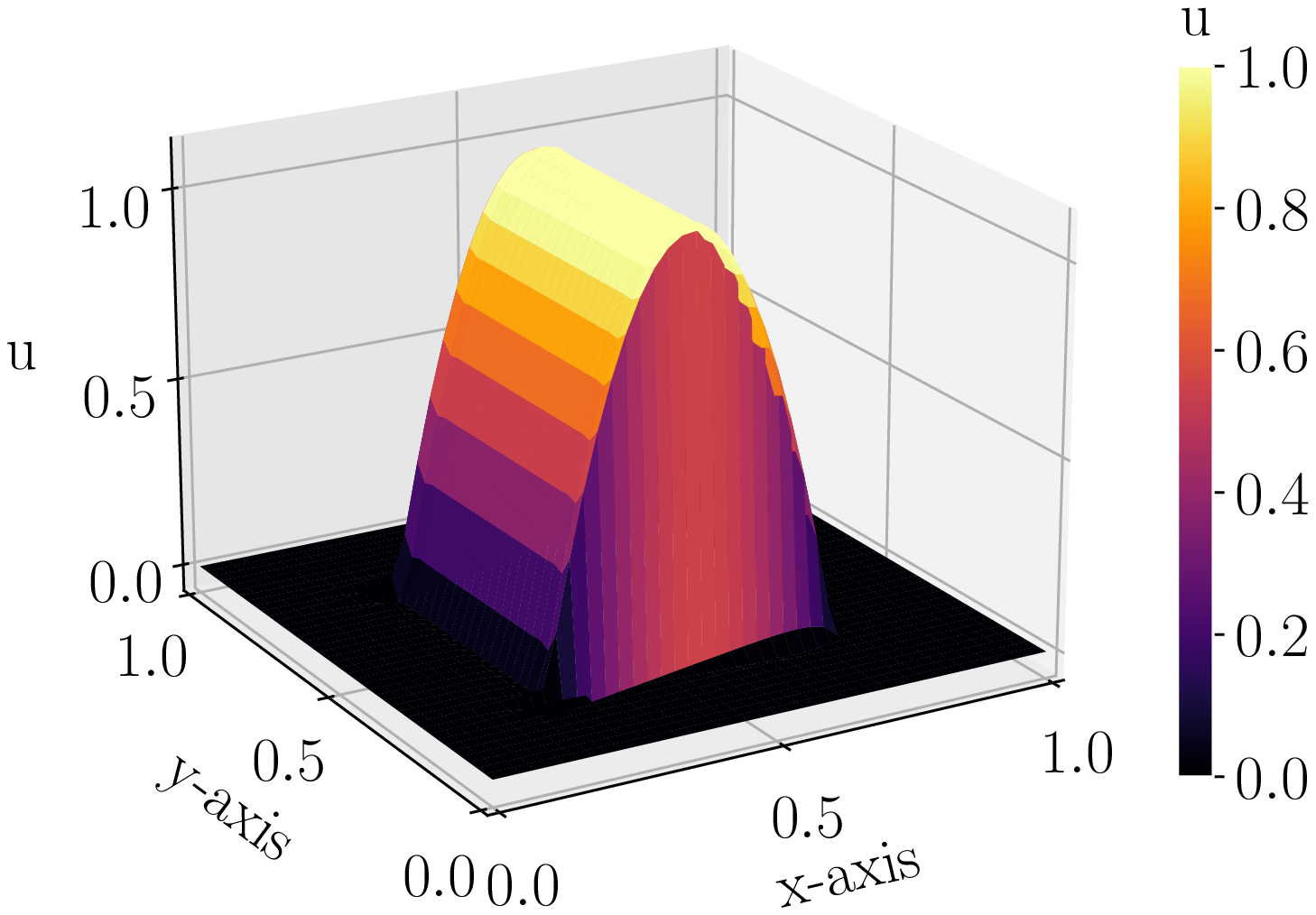}}
		\end{center}
    \caption{\label{Fig:std4} Solution with Standard $\tau_{std}$}
	\end{subfigure}%
	\newline
	\begin{subfigure}[t]{0.5 \textwidth}
		\begin{center}{\includegraphics[width=\textwidth]{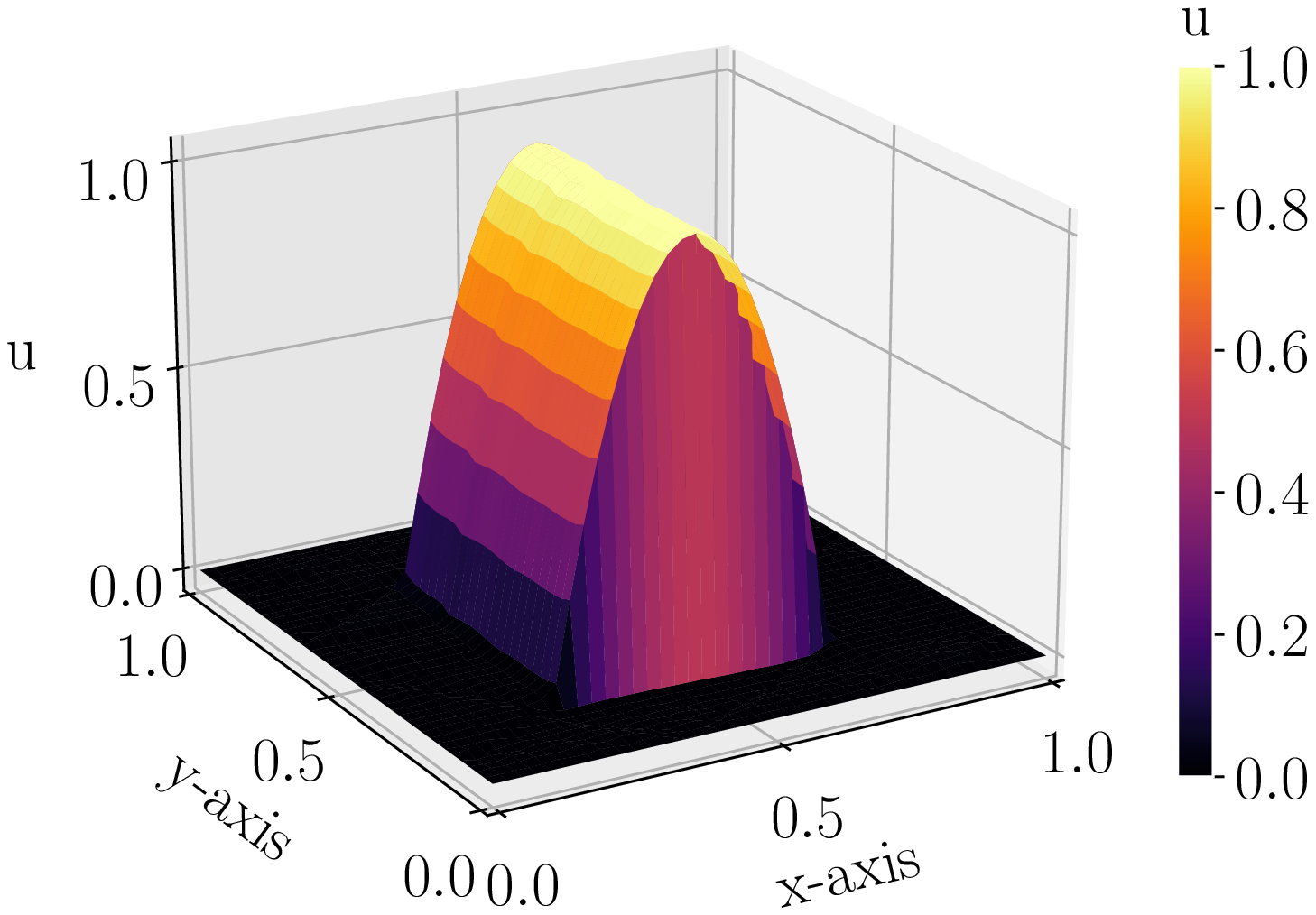}}
		\end{center}
    \caption{\label{Fig:sol4} Solution from AiStab-FEM~($u_h$)}
\end{subfigure}%
\begin{subfigure}[t]{0.5 \textwidth}
\begin{center}{\includegraphics[width=\textwidth]{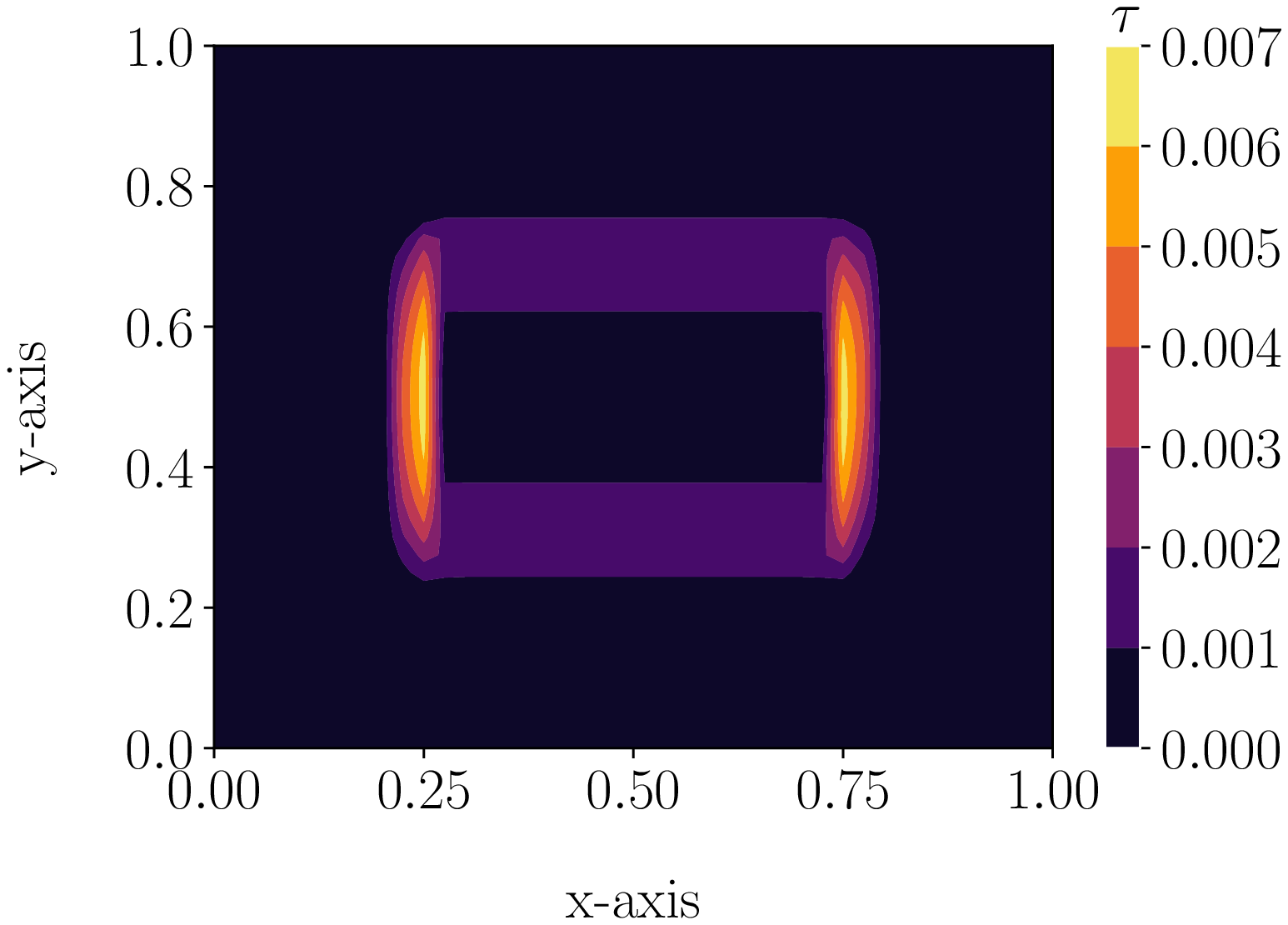}}
\end{center}
\caption{\label{Fig:tau4} $\tau_K$ predicted from AiStab-FEM}
\end{subfigure}%
\end{center}
\caption{\label{fig:Example4} Example 4}
\end{figure}
This example differs from example 1 in the source function $f$. It has been used in~\cite{articlepetr}.
It has two interior characteristic layers in the convection direction between the spatial points $(0.25,~0.25)$ and $(0.25,~0.75)$. The P\'eclet number is $1.77e+06$. We have solved this example for different mesh refinement levels to check the generalizability of AiStab-FEM. The numerical solution and heat map of $\tau_K$ corresponding to $h = \sqrt{2}/40$ are shown in figure \ref{fig:Example4}. The grid-convergence study has been presented in the next section. For this example, the solution obtained from AiStab-FEM is closest to the exact solution and has the least spurious oscillations. Further, AiStab-FEM performs better than other techniques on all the error metrics. AiStab-FEM and Standard $\tau_{std}$ capture the interior layers very well.
\subsection{Error Metrics}
We have compared the performance of VarNet, standard technique and Standard $\tau_{std}$ normalized with $||\nabla u_{Galerkin}||$ technique with the AiStab-FEM for all four benchmark examples in terms of $L^2$-error, $H^1$-semi-norm error, $L^\infty$-error and relative $l^2$-error.
PINN's performance is reported in relative $l^2$ error as other metrics are not computable on PINNs solution. The tables below show different error metrics for numerical solutions generated by considered numerical techniques.

\begin{table}[!htp]\centering
\caption{$L^2$-Error}\label{tab:l2error}
\begin{tabular}{lccccc}\toprule 
& &  \multicolumn{2}{c} {Examples}&\\\toprule 
Techniques
&1	&2	&3	&4	\\\midrule
Standard $\tau_{std}$	&1.32e-5&6.77e-6&1.41e-5&3.63e-6\\
Standard $\tau_{std}$ with $||\nabla u_{Galerkin}||$		&4.04e-5	&1.63e-5	&1.76e-5	&3.85e-6	\\
VarNet &1.70e-4&2.37e-4&3.00e-4&2.80e-4\\
AiStab-FEM	&6.05e-6	&5.04e-6	&1.20e-5	&3.63e-6	\\
\bottomrule
\end{tabular}
\end{table}

\begin{table}[!htp]\centering
\caption{Relative $l^2$-Error}\label{tab:rel2error}
\begin{tabular}{lrrrrr}
\toprule
& &  \multicolumn{2}{c} {Examples}&\\\toprule 
Techniques&1	&2	&3	&4	\\\midrule
PINN	&4.85e-1	&4.85e+1	&8.69e-1	&1.01e+0	\\
Standard $\tau_{std}$	&1.17e-1	&1.36e-1	&8.02e-2	&4.63e-2	\\
Standard $\tau_{std}$ with $||\nabla u_{Galerkin}||$	&3.35e-1	&3.06e-1	&9.68e-2	&4.90e-2	\\
VarNet &5.17e-1	&
1.62e+0	&
5.39e-1	&
1.25e+0	\\
AiStab-FEM	&6.17e-2	&9.73e-2	&6.94e-2	&4.60e-2	\\
\bottomrule
\end{tabular}
\end{table}

\begin{table}[!htp]\centering
\caption{$H^1$-Error}\label{tab:h1error}
\begin{tabular}{lrrrrr}\toprule
& &  \multicolumn{2}{c} {Examples}&\\\toprule 
Techniques&1	&2	&3	&4	\\\midrule
Standard $\tau_{std}$ & 1.30e-3	&6.74e-4	&1.43e-3	&3.28e-4	\\
Standard $\tau_{std}$ with $||\nabla u_{Galerkin}||$	&3.65e-3	&1.51e-3	&1.72e-3	&3.69e-4	\\
VarNet &1.80e-03	&
1.87e-3	&
2.26e-3	&
2.62e-3	\\
AiStab-FEM	&5.98e-4	&4.80e-4	&1.23e-3	&3.29e-4	\\
\bottomrule
\end{tabular}
\end{table}

\begin{table}[!htp]\centering
\caption{$L^\infty$-Error}\label{tab:linfinityerror}
\begin{tabular}{lrrrrr}
\toprule
& &  \multicolumn{2}{c} {Examples}&\\\toprule 
Techniques&1	&2	&3	&4	\\\midrule
Standard $\tau_{std}$ &9.07e-5	&7.29e-5	&1.13e-4	&7.60e-6	\\
Standard $\tau_{std}$ with $||\nabla u_{Galerkin}||$	&1.32e-4	&8.70e-5	&1.47e-4	&1.06e-5	\\
VarNet &3.27e-4	&3.55e-4	&3.53e-4	&3.54e-4	\\
AiStab-FEM	&3.99e-5	&4.05e-5	&9.08e-5	&7.35e-6	\\
\bottomrule
\end{tabular}
\end{table}
\vskip 20pt
\subsection{Mesh refinement analysis}
Studying grid convergence is essential for ensuring any mesh-based solver's numerical stability and accuracy. Thus we show the numerical errors obtained on refining the mesh subsequently, for example 4. We compute the order of error convergence for the $L^2$-Error and obtain the optimal order of convergence for $P_2$ finite element for $u$,   as shown in Table \ref{tab:mesh_refinement}. 
\begin{table}[H]\centering
\caption{Errors in the solution}\label{tab: }
\begin{tabular}{lcrrrrrr}\toprule
$N\_cells$ &h &Residual &$L^2$-Error &Relative &$H^1$-Error &$L^{\infty}$-Error \\
 & & & &$l^2$-Error & & \\\midrule
10 &1.41e-1 &8.71e-1 &4.75e-4 &3.83e-1 &1.06e-2 &1.04e-3 \\
20 &7.07e-2 &5.12e-1 &2.93e-5 &9.28e-2 &1.34e-3 &5.92e-5 \\
40 &3.54e-2 &2.88e-1 &3.63e-6 &4.60e-2 &3.29e-4 &7.35e-6 \\
80 &1.77e-2 &1.76e-1 &4.51e-7 &2.29e-2 &8.14e-5 &9.20e-7 \\
\bottomrule
\end{tabular}
\end{table}

\begin{table}[!htp]\centering
\caption{Mesh refinement analysis}
\begin{tabular}{lcrrr}\toprule
$N\_cells$ &h &$L^2$-Error &Order \\\midrule
10 &1.41e-1 &4.75e-4 & \\
20 &7.07e-2 &2.93e-5 &4.02 \\
40 &3.54e-2 &3.63e-6 &3.01 \\
80 &1.77e-2 &4.51e-7 &3.01 \\
\bottomrule
\end{tabular}
\label{tab:mesh_refinement}
\end{table}

\section{Conclusion}
This paper proposes a neural network-based technique to predict an optimal stabilization parameter for the SUPG method to solve SPPDEs. 
In this method, the equation coefficients and the mesh size are taken as input features 
for the neural network. 
The predicted $\hat{\tau}$ is normalized by $||\nabla u_{Galerkin}||$.
The variational stabilized form of the equation is used to calculate the cost for the 
back-propagation algorithm. 

The proposed technique outperforms the existing Neural Network-based PDE solvers such as PINNs and VarNet for application to SPPDEs. 
It also performs better than the standard technique for benchmark problems. 
The proposed technique achieves the optimal order of convergence when the computational grid is refined. 

\bibliographystyle{elsarticle-num} 
\bibliography{main}

\end{document}